\documentclass[reqno,11pt]{amsart}

\usepackage{amscd,amssymb,comment,epic,eepic,euscript}
\usepackage{graphicx}
\usepackage{enumerate}
\usepackage[initials,alphabetic]{amsrefs}
\usepackage{array}
\usepackage{longtable}
\usepackage{algorithm}
\usepackage[noend]{algpseudocode}
\usepackage[table,x11names]{xcolor}

\usepackage{listings}
\input{listings-python.prf}
\usepackage{pgf,tikz}

\usepackage{mathtools}

\usepackage{enumitem,color,titletoc}%

\definecolor{dkgreen}{rgb}{0,0.6,0}
\definecolor{gray}{rgb}{0.5,0.5,0.5}
\definecolor{mauve}{rgb}{0.58,0,0.82}
\definecolor{purple2}{RGB}{153,0,153} %
\definecolor{green2}{RGB}{0,153,0} %
\definecolor{ogreen}{RGB}{85,107,47} %
\definecolor{dred}{RGB}{139,0,0}

\def\bes{\begin{equation*} }
\def\ees{\end{equation*} }

\usepackage{xspace}

\newcommand{\esorou}{equisigned sum of roots of unity\xspace}

\usepackage{multicol}

\renewcommand{\thetable}{\arabic{table}}

\theoremstyle{plain}
\newtheorem{thm}{Theorem}[section]
\newtheorem{lemma}[thm]{Lemma} %
\newtheorem{prop}[thm]{Proposition}
\newtheorem{conj}[thm]{Conjecture}
\newtheorem{algo}[thm]{Algorithm} %

\theoremstyle{remark}

\newtheorem{remark}[thm]{Remark}

\theoremstyle{definition}
\newtheorem{defi}[thm]{Definition}

\newcount\theTime
\newcount\theHour
\newcount\theMinute
\newcount\theMinuteTens
\newcount\theScratch
\theTime=\number\time
\theHour=\theTime
\divide\theHour by 60
\theScratch=\theHour
\multiply\theScratch by 60
\theMinute=\theTime
\advance\theMinute by -\theScratch
\theMinuteTens=\theMinute
\divide\theMinuteTens by 10
\theScratch=\theMinuteTens
\multiply\theScratch by 10
\advance\theMinute by -\theScratch
\def\timeHHMM{{\number\theHour:\number\theMinuteTens\number\theMinute}}

\def\today{{\number\day\space
 \ifcase\month\or
  January\or February\or March\or April\or May\or June\or
  July\or August\or September\or October\or November\or December\fi
 \space\number\year}}

\def\timeanddate{{\timeHHMM\space o'clock, \today}}

\newcommand{\Z}{{\mathbb Z}}
\newcommand{\Q}{{\mathbb Q}}

\newcommand{\N}{{\mathbb N}}

\newcommand\val{\operatorname{val}}

\def\bes{\begin{equation*} }
\def\ees{\end{equation*} }

\def\beq{ \begin{equation} }
\def\eeq{ \end{equation} }

\def\bep{\begin{proof}}
\def\eep{\end{proof}}

\def\ben{\begin{enumerate}}
\def\een{\end{enumerate}}

\def\bet{\begin{theorem}}
\def\eet{\end{theorem}}

\def\bel{\begin{lemma}}
\def\eel{\end{lemma}}

 \usepackage[normalem]{ulem}

\usepackage[colorlinks,linkcolor=dred,citecolor=ogreen,urlcolor=gray,hypertexnames=true]{hyperref}

\begin{document}

\title[Vanishing sums of roots of unity]{Classifying minimal vanishing\\ sums of roots of unity}

\author[Christie]{Louis Christie} 
 \address{L.C., Department of Mathematics, 
 The University of Auckland, Private Bag 92019, Auckland 1142, New Zealand}
  \email{louis@christie.kiwi.nz}

\author[Dykema]{Kenneth J. Dykema$^{*}$}
 \address{K.D., Department of Mathematics, Texas A\&M University,
  College Station, TX 77843-3368, USA}
 \email{kdykema@math.tamu.edu}
 \thanks{${}^{*}$Research supported in part by NSF grant DMS--{1800335 and by a grant from the Simons Foundation/SFARI (524187, K.D.)}}
\author[Klep]{Igor Klep${}^{\star}$}
 \address{I.K., Faculty of Mathematics and Physics,
 University of Ljubljana, Jadranska 21, 1000 Ljubljana, Slovenia \& 
Faculty of Mathematics, Natural Sciences
and Information Technologies, 
University of Primorska, 
Glagoljaška 8, 6000 Koper, Slovenia \&
Institute of Mathematics, Physics and Mechanics,
Jadranska 19, 1000 Ljubljana, Slovenia}
 \email{igor.klep@fmf.uni-lj.si}
 \thanks{${}^\star$Supported by the Slovenian Research and Innovation Agency 
program P1-0222 and
 grants J1-50002, N1-0217, J1-60011, J1-50001, J1-3004 and J1-60025.
Partially supported by the Marsden Fund Council of the Royal Society of New Zealand.}
\date{\timeanddate}

\setcounter{tocdepth}{2}
\contentsmargin{2.55em} 
\dottedcontents{section}[3.8em]{}{2.3em}{.4pc} 
\dottedcontents{subsection}[6.1em]{}{3.2em}{.4pc}

\makeatletter
\newcommand{\mycontentsbox}{%
{%
\linespread{1.2}
\small\tableofcontents}}
\def\enddoc@text{\ifx\@empty\@translators \else\@settranslators\fi
\ifx\@empty\addresses \else\@setaddresses\fi
\newpage\mycontentsbox}
\makeatother

 \begin{abstract}
A vanishing sum of roots of unity 
is called
minimal if 
no proper, nonempty sub-sum of it vanishes.
This paper classifies all minimal vanishing sums of roots of unity of weight $\leq 16$ by hand,
{thereby uncovering new phenomena beyond the earlier 1998 classification of Poonen and Rubinstein (SIAM J. Discrete Math.) that went up to weight $12$.}
{The paper also develops an algorithm to explore higher weights up to $21$, yielding a conjectural extension of the classification.}
\looseness=-1
\end{abstract}

\maketitle

\section{Introduction}

A {\em root of unity} is a complex number $\omega$, some positive integer power of which is equal to $1$.
This paper studies vanishing sums of roots of unity, i.e., integer combinations 
\begin{equation}\label{eq:sorouIntro}
\sum_{i=1}^k a_i \omega_i,
\end{equation} where $a_i\in\Z$
and $\omega_i$ are
{distinct}
roots of unity,
{so that the value of the sum}
is $0$.
If each proper sub-sum of such a sum is nonzero, we call \eqref{eq:sorouIntro} {\em minimal}. 
The \emph{weight} of 
an integer combination of the form
\eqref{eq:sorouIntro} is
defined to be 
$\ell_1$-norm of the coefficient vector, i.e., 
$\sum_i |a_i|$.
The {\em height} is the $\ell_\infty$-norm,
{namely,}
$\max_i |a_i|$.

Vanishing sums of roots of unity naturally arise 
 in a number of areas in algebra \cites{Sch64,Ste08,AS12}, geometry and combinatorics \cite{PR98}, and number theory \cites{CJ76,DZ00,Len79,Man65};
 we refer the reader to \cite{LL00} for more details and further references.
 However, they are ubiquitous across many sciences. For instance, they {have} appeared
 in operator algebra \cite{JMS14}, nuclear and particle physics \cite{FG14} and approximation theory \cite{CW11}. 
Mann \cite{Man65} classifies all minimal
vanishing sums of roots of unity of weight
$\leq7$, Conway and Jones \cite{CJ76}
extend this to weight $\leq9$, and
Poonen and Rubinstein \cite{PR98}
classify the ways in which 12 roots of unity can sum to zero.
They introduced a notion of {\em type} of a minimal vanishing sum of roots of unity.

{
Our detailed investigations of vanishing sums of roots of unity grew out of an attempt to prove the Kaplansky-Lvov conjecture in degree 3.
We needed a catalog of the minimal vanishing sums of roots of unity.
It {was important}, for our purpose, to understand their types and their {parities}, namely, the number of roots of unity
of odd and even order in the vanishing sum.
Thus, in this paper, we extend the classification of Poonen and Rubinstein~\cite{PR98} by extending the definition of type and %
{implementing an algorithm that produces, up to weight $21$, a list of types and parities of minimal vanishing sums of roots of unity. This list agrees with known classifications up to weight $16$ and suggests a conjectural classification up to weight $21$.}
}

\subsection{Contributions and reader's guide}
Section \ref{sec:not} fixes notation and terminology, and presents a few of the preliminaries needed later in the paper. 
Our main result, Theorem \ref{thm:mvs} in Section \ref{sec:hand}, 
classifies all such
{types and parities for}
all minimal vanishing sums of roots of unity of
weight $\leq16$. 
There are 76 types in this classification,
and the number of types increases rapidly with weight.
The proof of Theorem~\ref{thm:mvs} does not rely on use of computers.
New phenomena arise when extending the classification beyond weight $12$, see
Proposition \ref{prop:hnuf}, Definition \ref{defi:type} and Subsection \ref{ssec:ex} for details.

Our computational exploration suggests that the lowest weight minimal vanishing sums of roots of unity of height $>1$ occur at weight $21$; see Conjecture \ref{conj:ht2} and the data in Table \ref{tab:mvs21}. The large number of cases to handle forces us to turn to computer assistance: in Section \ref{sec:comp1} we present an algorithm that, in principle, finds all minimal vanishing sums of a given weight, 
describe its implementation,
and report the resulting list up to weight $21$ (Appendix \ref{app:all}). Since we do not formally verify the implementation, the conclusions of Section \ref{sec:comp1} should be viewed as conjectural.

\section{Notation and Preliminaries}\label{sec:not}

For $n\in\N$, we write $\nu_n=e^{2\pi i/n}$ for the standard primitive $n$-th root of unity.
The {\em order} of a root of unity $\omega$ is the least positive integer $d$ such that $\omega^d=1$.

A {\em sum of roots of unity} (also written {\em sorou}) is an
unordered, finite, nonempty list $h=(\omega_1,\omega_2,\ldots,\omega_n)$ of roots of unity.\footnote{The plural of sorou is {\em sums of roots of unity}, which is also written as {\em sorou}.}
The {\em terms} of $h$ are the roots of unity $\omega_1,\ldots,\omega_n$.
The {\em value} of $h$ is the complex number
\[
\val(h)=\sum_{j=1}^n\omega_j
\]
and the {\em weight} of $h$ is the integer $n$.
The {\em rotation} of $h$ by the root of unity $z$ is $zh=(z\omega_1,\ldots,z\omega_n)$.
{We write $h\sim zh$ and note that $\sim$ is an equivalence relation.}
The {\em order} of $h$ is the least common multiple of the orders of the terms $\omega_j$ of $h$,
and the {\em relative order} of $h$ is the least common multiple of the orders of all the ratios $\omega_i/\omega_j$ of all the terms of $h$.
Note that a sorou with relative order $d$ can always be rotated to obtain a sorou of order $d$.

As is traditional, we will usually write $\omega_1+\cdots+\omega_n$ for the sorou, instead of writing it as a list.
With this convention, the sorou $(1,-1)$ is represented as $1+(-1)$, which is not the same as the complex number $0$, though its value is,
of course $0$.
As is natural with this notation, for sorou $f$ and $g$ we let $f+g$ denote the sorou that, technically,
is the concatenation of the lists of terms of $f$ and of $g$.
The {\em multiplicity} of a root of unity $z$ in a sorou $h$ as above is the number of times it appears in $h$,
namely, the number of
$j\in\{1,\ldots,n\}$ such that $\omega_j=z$.  We may write a sorou as a sum employing multiplicities.
Thus, $2\nu_3+\nu_3^2+3(-1)$ represents the sorou $h=(\nu_3,\nu_3,\nu_3^2,-1,-1,-1)$.
The {\em height} of $h$ is the maximum multiplicity of any of the $\omega_j$ that appears in $h$.
We also employ natural conventions regarding minus signs, so that the same sorou
{$h$}
may
{also}
be written $2\nu_3+\nu_3^2-3$,
and $f-g$ means the sorou $f+(-1)g$.

A {\em subsorou} of a sorou $h=(\omega_1,\ldots,\omega_n)$ is a sorou $g$ 
that can be written as
$g=(\omega_{k(1)},\omega_{k(2)},\ldots,\omega_{k(m)})$ for some $1\le k(1)<\cdots< k(m)\le n$,
for some $m\in\{0,\ldots,n\}$.
If $m=0$, then we say $g$ is the {\em empty subsorou}, and if $m<n$, i.e., $g\ne h$, then 
we say $g$ is a {\em proper subsorou}.
{We write $g\prec h$ when $g$ is a subsorou of $h$.}

A sorou is said to {\em vanish} if its value is $0$.
A vanishing sorou is said to be {\em minimal vanishing} if no proper, nonempty subsorou of it is vanishing.
It is clear that every vanishing sorou can be written as a sum (or equivalently, in the list notation, as a concatenation)
of minimal vanishing sorou, though not necessarily in a unique way.

As mentioned in the introduction, vanishing sorou have been studied by many authors.
We will classify the minimal vanishing sorou of weight $\le16$,
{and provide computational evidence towards a classification up to weight $21$.}
This continues and extends the work~\cite{PR98} of Poonen and Rubinstein, who classified the minimal vanishing sorou of weight
no more than $12$.

The following lemma is a consequence of \cite[Theorem 1]{Man65}.
\begin{lemma}
The relative order of a minimal vanishing sorou $h$ is a product $p_1p_2\cdots p_s$ of distinct primes $p_1<p_2<\cdots<p_s$.
\end{lemma}

\begin{defi}
In the setting of the above lemma, we call $p_s$ the {\em top prime} of $h$.
\end{defi}

\begin{prop}\label{prop:hnuf}
Let $h$ be a sorou whose order is a product $p_1p_2\ldots p_s$ of distinct primes $p_1<p_2<\cdots<p_s$.
Let $p=p_s$ be the top prime.
Then, after replacing $h$ by a rotation, if necessary,
we have
\begin{equation}\label{eq:hnuf}
h=\sum_{j=0}^{p-1}\nu_p^jf_j,
\end{equation}
for some sorou $f_0,\ldots,f_{p-1}$, each term of which has order dividing $p_1p_2\cdots p_{s-1}$.
Then $h$ vanishes if and only if
\begin{equation}\label{eq:valfs}
\val(f_0)=\val(f_1)=\cdots=\val(f_{p-1}).
\end{equation}

Suppose that $h$ vanishes.
Then $h$ is minimal vanishing if and only if the following hold:
\begin{enumerate}[label={\rm(\roman*)}]
\item $\val(f_0)\ne0$,
\item for no $j$ does $f_j$ have a vanishing proper, nonempty, subsorou,
\item there is no complex number $z$ such that for all $j$, $f_j$ has a proper, nonempty subsorou
with value $z$.
\end{enumerate}
\end{prop}

\begin{proof}
The first statement follows from the fact that 
$\Q(\nu_{p},\nu_{p_{s-1}}\cdots\nu_{p_1})/\Q(\nu_{p_{s-1}}\cdots\nu_{p_1})$
is a field extension of degree $p$,
and the minimal polynomial for $\nu_p$ is
$\Phi_{p}(x)=\sum_{j=0}^{p-1}x^j$ (cf.~\cite[Proposition 2.4]{Wah97}).
Then the minimality characterization follows easily.
\end{proof}

Note that equation~\eqref{eq:valfs} is equivalent to the fact that for all $i,j\in\{0,\ldots,p-1\}$, the sorou $f_i-f_j$ vanishes.

From Proposition~\ref{prop:hnuf}, we can deduce the well known fact (cf.~Definition \ref{defi:type}\ref{it:Rp} below)
that if $h$ is a minimal vanishing sorou of prime order $p$, then
up to rotation $h$ equals
\[
1+\nu_p+\nu_p^2\cdots+\nu_p^{p-1}.
\]
Following~\cite{PR98}, we say that such a vanishing sorou has type $R_p$.

\subsection{Types and parities}

\begin{defi}\label{defi:type}
We recursively associate to an arbitrary vanishing sorou $h$
(or, more correctly, to the equivalence class of $h$ under rotation)
{its} {\em type} {(or types)} as follows.
\begin{enumerate}[label={\rm(\alph*)}]
\item
If $h$ is a sum of minimal vanishing sorou $v_1,\ldots,v_k$ having types $T_1,\ldots,T_k$,
respectively, then we say $h$ has type $T_1\oplus T_2\oplus\cdots\oplus T_k$.
\item\label{it:subsidsorou}
A minimal vanishing sorou $h$, after rotation, can always be written as in~\eqref{eq:hnuf}
with $w(f_0)\le w(f_j)$ for all $j\le p-1$ and and with $1\prec f_0$.
We call $f_0,\ldots,f_{p-1}$ the {\em subsidiary sorou} of $h$.
\item
The weights $w(f_0),\ldots,w(f_{p-1})$, when arranged in increasing order,
form the {\em subsidiary weight partition} of $h$, often called just the {\em weight partition} of $h$.
\item
We call $f_0$ a {\em smallest weight subsidiary sorou} of $h$.
Note that $f_0$ is not, in general, unique,
{even up to rotation,}
since we allow rotations of $h$, but the weight of $f_0$ is unique,
and is called the {\em smallest subsidiary weight} of $h$.
\item\label{it:Rp} 
If $f_j=f_0$ for all $j$, then $f_0$ can have only one term, namely, $1$, (for otherwise, by Proposition~\ref{prop:hnuf}, $h$
would fail to be minimal vanishing); in this case, as noted already, we say $f$ has type $R_p$.
\item
For minimal vanishing $h$ written (after rotation) as in~\ref{it:subsidsorou}, let $J$
be the set of all $j\in\{1,\ldots,p-1\}$ for which $f_j\ne f_0$ and suppose $J=\{j(1),j(2),\ldots,j({n})\}$ is nonempty.
Suppose the vanishing sorou $f_0-f_{j(i)}$ has type $T_i$.
Then we say that $h$ has type 
\begin{equation}\label{eq:type}
(R_p:f_0:T_1,T_2,\ldots,T_n),
\end{equation}
where the ordering of $T_1,\ldots,T_n$ is unimportant.
We call $T_1,\ldots,T_n$ the {\em subsidiary types} of $h$.
\item In the case $f_0=1$, we may omit the ``$f_0:$'', writing instead
\[
(R_p:T_1,T_2,\ldots,T_n).
\]
Furthermore, whenever some of the subsidiary types $T_i$ are repeated, the multiplicity may be indicated with an integer.
For example, $(R_7:3R_3,2R_5)$ instead of $(R_7:R_3,R_3,R_3,R_5,R_5)$.
\end{enumerate}
\end{defi}

\begin{remark}\label{rem:type}
\mbox{}\par
\begin{enumerate}[label={\rm(\alph*)}]
\item The above definition builds upon the definition used in~\cite{PR98} (which was confined to the case $f_0=1$).
\item A given vanishing sorou may have more than one type.
For example $1+\nu_3+\nu_3^2-1-\nu_3-\nu_3^2$ has type $R_3\oplus R_3$ and type $R_2\oplus R_2\oplus R_2$.
\item In some cases, for example $R_p$ or $(R_p:R_q)$, the type specifies the sorou uniquely up to rotation.
But not in other cases, for example $(R_5:2R_3)$.
\item The types of minimal vanishing sorou are distinct from the types of non-minimal vanishing sorou.
\item If a minimal vanishing sorou has type $(R_p:T_1,\ldots,T_n)$, namely, with smallest subsidiary weight $1$,
then all of the types $T_i$ must be of minimal vanishing sorou.
\item\label{it:Tj} More generally, if the type is $(R_p:f_0:T_1,\ldots,T_n)$, then each $T_j$ decomposes as a sum of at most $w(f_0)$ minimal vanishing
sorou (for otherwise, $f_j$ would itself have a proper, nonempty, vanishing subsorou).
\end{enumerate}
\end{remark}

We also define the \textit{parity} of a sorou as the pair of integers counting the number of positive and negative signs in the sorou. 
Note that because we allow rotations, the order of the pair is irrelevant.

\subsection{A few examples}\label{ssec:ex}

Now that we have introduced basic terminology, let us illustrate it with a few examples. 
Consider the sorou $h_1 = 1 + \nu_5 + \nu_5^2 + \nu_5^3 + \nu_5^4$ of type $R_5$ and the sorou $h_2 = 1 + \nu_3 + \nu_3^2$ of type $R_3$. 
{Then $h_1'=\nu_5 + \nu_5^2 + \nu_5^3 + \nu_5^4$ and $h_2'=\nu_3 + \nu_3^2$ are both sorou having valuation $-1$, and subtracting them yields
the minimal vanishing sorou}
\[ h = {h_1' - h_2'} = \nu_5 + \nu_5^2 + \nu_5^3 + \nu_5^4 - \nu_3 - \nu_3^2, \]
{which can be rewritten}
\[ h = \nu_5 + \nu_5^2 + \nu_5^3 + \nu_5^4 + \nu_6 + \nu_6^5 \sim \sum_{j = 0}^4 f_j \nu_5^j \]
where $f_0 = f_1 = f_2 = f_3 = 1$ and $f_4 = \nu_6 + \nu_6^5$. Thus $f_0 - f_4 = 1 - \nu_6 - \nu_6^5 = 1 + \nu_3 + \nu_3^2$ and so we have a minimal vanishing sorou, with type $(R_5 : R_3)$.
{In a similar manner, we obtain
\[
g=\nu_5^2+\nu_5^3+\nu_5^4-\nu_3-\nu_3^2-\nu_5\nu_3-\nu_5\nu_3^2
\]
of type}
$(R_5:2R_3)$, etc.
{These are pictured below.}
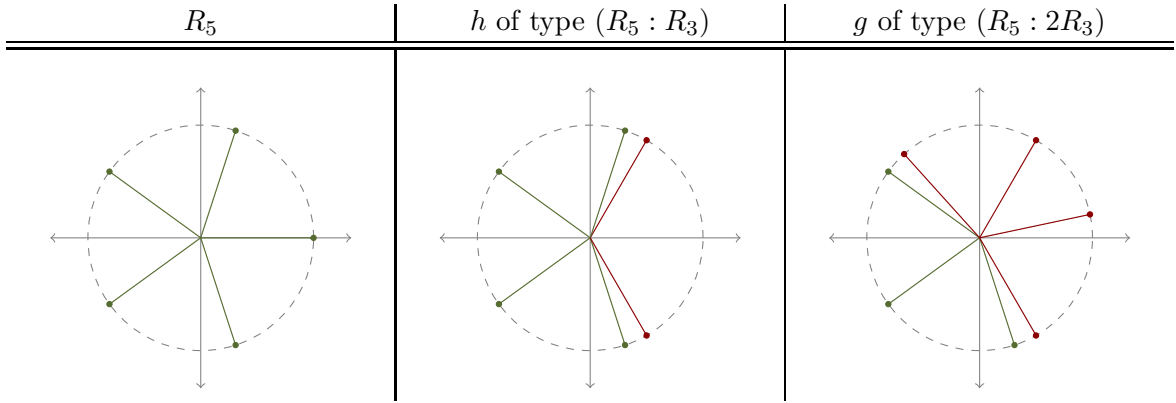
\begin{center}
\begin{tabular}{c|c|c}
$R_5$ & {$h$ of type} $(R_5 : R_3)$ & {$g$ of type} $(R_5 : 2 R_3)$  \\ \hline \hline
$\phantom{0000000000000000000000000}$ & $\phantom{0000000000000000000000000}$ & $\phantom{0000000000000000000000000}$ \\
\begin{tikzpicture}[scale = 0.5]
\draw[<->, gray] (-4,0) -- (4,0);
\draw[<->, gray] (0,-4) -- (0,4);
\draw[dashed, thin, gray] (0,0) circle (3);
\draw[fill = ogreen, color=ogreen] (0,0) -- (0 : 3) circle (0.7mm);
\draw[fill = ogreen, color=ogreen] (0,0) -- (72 : 3) circle (0.7mm);
\draw[fill = ogreen, color=ogreen] (0,0) -- (144 : 3) circle (0.7mm);
\draw[fill = ogreen, color=ogreen] (0,0) -- (216 : 3) circle (0.7mm);
\draw[fill = ogreen, color=ogreen] (0,0) -- (288 : 3) circle (0.7mm);
\end{tikzpicture} & \begin{tikzpicture}[scale = 0.5]
\draw[<->, gray] (-4,0) -- (4,0);
\draw[<->, gray] (0,-4) -- (0,4);
\draw[dashed, thin, gray] (0,0) circle (3);
\draw[fill = ogreen, color=ogreen] (0,0) -- (144 : 3) circle (0.7mm);
\draw[fill = ogreen, color=ogreen] (0,0) -- (216 : 3) circle (0.7mm);
\draw[fill = ogreen, color=ogreen] (0,0) -- (288 : 3) circle (0.7mm);
\draw[fill = dred, color=dred] (0,0) -- (300 : 3) circle (0.7mm);
\draw[fill = dred, color=dred] (0,0) -- (60 : 3) circle (0.7mm);
\draw[fill = ogreen, color=ogreen] (0,0) -- (72 : 3) circle (0.7mm);
\end{tikzpicture} & \begin{tikzpicture}[scale = 0.5]
\draw[<->, gray] (-4,0) -- (4,0);
\draw[<->, gray] (0,-4) -- (0,4);
\draw[dashed, thin, gray] (0,0) circle (3);
\draw[fill = ogreen, color=ogreen] (0,0) -- (144 : 3) circle (0.7mm);
\draw[fill = ogreen, color=ogreen] (0,0) -- (216 : 3) circle (0.7mm);
\draw[fill = ogreen, color=ogreen] (0,0) -- (288 : 3) circle (0.7mm);
\draw[fill = dred, color=dred] (0,0) -- (300 : 3) circle (0.7mm);
\draw[fill = dred, color=dred] (0,0) -- (60 : 3) circle (0.7mm);
\draw[fill = dred, color=dred] (0,0) -- (12 : 3) circle (0.7mm);
\draw[fill = dred, color=dred] (0,0) -- (132 : 3) circle (0.7mm);
\end{tikzpicture}
\end{tabular}
\end{center}
{Now consider the minimal vanishing sorou $\sum_{j=o}^6 \nu_7^j$. Multiplying by $1+\nu_5$ yields another vanishing, but not minimal sorou, but then subtracting $h_1$ gives:
\[
\Big(\sum_{i=0}^6 \nu_7^i\Big)(1+\nu_5)-\Big(\sum_{i=0}^4\nu_5^j\Big) \sim \sum_{j=0}^6 f_j \nu_7^j,
\]
where
\[
f_0 = f_1 = \cdots = f_5 = 1+ \nu_5 \qquad \text{and} \qquad f_6 = \nu_{10} + \nu_{10}^3 + \nu_{10}^9,
\]
a minimal vanishing sorou of type} $(R_7 : 1 + \nu_5 : R_5)$, with $f_0\succ1$.
For types with $f_0\succ1$ we can have non-minimal subtypes
(see Remark \ref{rem:type}\ref{it:Tj}) as depicted below.

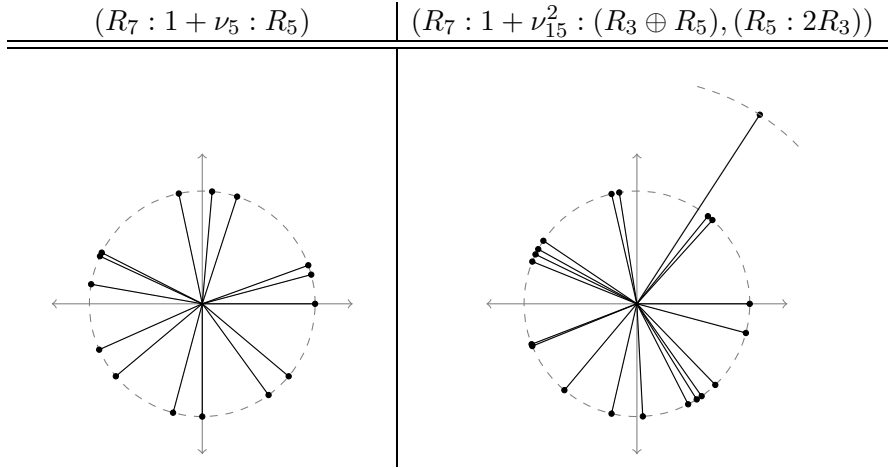
\begin{center}
\begin{tabular}{c|c}
$(R_7 : 1 + \nu_5 : R_5)$ & $(R_7 : 1 + \nu_{15}^2 : (R_3 \oplus R_5), (R_5 : 2R_3))$ \\ \hline \hline
$\phantom{0000000000000000000000000}$ & $\phantom{0000000000000000000000000}$ \\

\begin{tikzpicture}[scale = 0.5]
\draw[<->, gray] (-4,0) -- (4,0);
\draw[<->, gray] (0,-4) -- (0,4);
\draw[dashed, thin, gray] (0,0) circle (3);
\draw[fill = black] (0,0) -- (0:3) circle (0.7mm);
\draw[fill = black] (0,0) -- (72:3) circle (0.7mm);
\draw[fill = black] (0,0) -- (102:3) circle (0.7mm);
\draw[fill = black] (0,0) -- (153:3) circle (0.7mm);
\draw[fill = black] (0,0) -- (204:3) circle (0.7mm);
\draw[fill = black] (0,0) -- (255:3) circle (0.7mm);
\draw[fill = black] (0,0) -- (306:3) circle (0.7mm);
\draw[fill = black] (0,0) -- (20:3) circle (0.7mm);
\draw[fill = black] (0,0) -- (170:3) circle (0.7mm);
\draw[fill = black] (0,0) -- (220:3) circle (0.7mm);
\draw[fill = black] (0,0) -- (270:3) circle (0.7mm);
\draw[fill = black] (0,0) -- (320:3) circle (0.7mm);
\draw[fill = black] (0,0) -- (15:3) circle (0.7mm);
\draw[fill = black] (0,0) -- (85:3) circle (0.7mm);
\draw[fill = black] (0,0) -- (155:3) circle (0.7mm);
\end{tikzpicture} &

\begin{tikzpicture}[scale=0.5]
\draw[<->, gray] (-4,0) -- (4,0);
\draw[<->, gray] (0,-4) -- (0,4);
\draw[dashed, thin, gray] (0,0) circle (3);
\draw[fill = black] (0,0) -- (0:3) circle (0.7mm);
\draw[fill = black] (0,0) -- (48:3) circle (0.7mm);
\draw[fill = black] (0,0) -- (51:3) circle (0.7mm);
\draw[fill = black] (0,0) -- (99:3) circle (0.7mm);
\draw[fill = black] (0,0) -- (103:3) circle (0.7mm);
\draw[fill = black] (0,0) -- (151:3) circle (0.7mm);
\draw[fill = black] (0,0) -- (154:3) circle (0.7mm);
\draw[fill = black] (0,0) -- (202:3) circle (0.7mm);
\draw[fill = black] (0,0) -- (146:3) circle (0.7mm);
\draw[fill = black] (0,0) -- (314:3) circle (0.7mm);
\draw[fill = black] (0,0) -- (158:3) circle (0.7mm);
\draw[fill = black] (0,0) -- (230:3) circle (0.7mm);
\draw[fill = black] (0,0) -- (302:3) circle (0.7mm);
\draw[fill = black] (0,0) -- (257:3) circle (0.7mm);
\draw[fill = black] (0,0) -- (305:3) circle (0.7mm);
\draw[fill = black] (0,0) -- (201:3) circle (0.7mm);
\draw[fill = black] (0,0) -- (273:3) circle (0.7mm);
\draw[fill = black] (0,0) -- (345:3) circle (0.7mm);
\draw[fill = black] (0,0) -- (57:6) circle (0.7mm);
\draw[fill = black] (0,0) -- (297:3) circle (0.7mm);
\begin{scope}
\clip (-1,2) rectangle (4.3,5.8);
\draw[dashed, thin, gray] (0,0) circle (6);
\end{scope}
\end{tikzpicture}

\end{tabular}
\end{center}
In fact, the minimal vanishing sorou
{
$h= \sum_{j=0}^6\nu_7^jf_j$ shown above
with
\begin{align*}
f_j&=1+\nu_3\nu_5^4\quad(j\in\{0,1,2,3,5\}) \\
f_4&=-\nu_3-\nu_5^4+\nu_3^2(\nu_5+\nu_5^2+\nu_5^3) \\
f_6&=-\nu_5-\nu_5^2-\nu_5^3-2\nu_5^4-\nu_3^2\nu_5^4
\end{align*}
of height $2$ and weight $21$ is of}
type $(R_7 : 1 + \nu_{15}^2 : (R_3 \oplus R_5), (R_5 : 2R_3))$.
{Our computations in Section \ref{sec:comp1} suggest that this is one of only five types of minimal vanishing sorou having height $>1$ and weight $\leq21$ (see Conjecture \ref{conj:ht2}).}

Verification that $h$ is minimal vanishing using the criteria of Proposition~\ref{prop:hnuf}
can be done ``by hand'', symbolically or numerically using both Mathematica and our python code.

\section{Classification of sorou of weight $\leq16$}\label{sec:hand}

The main result of this section 
is Theorem \ref{thm:mvs} that 
classifies minimal vanishing 
sorou of weight $\leq16$.
Theorem \ref{thm:2pq}
presents a simple characterization
of minimal vanishing sorou whose
relative order divides $2pq$.

The following is Lemma~3.3 of~\cite{PR98} (and an immediate consequence of it).
\begin{lemma}\label{lem:2p}
Let $p$ be an odd prime.
The only minimal vanishing sorou of relative order dividing $2p$ are of types $R_2$ and $R_p$.
Thus, the only vanishing sorou of relative order dividing $2p$ are of types $R_2^{\oplus m}\oplus R_p^{\oplus n}$ for $n,m\ge0$.
\end{lemma}

\begin{thm}\label{thm:2pq}
Let $p<q$ be odd primes.
Suppose $h$ is a minimal vanishing sorou whose relative order divides $2pq$.
Then either $h$ is of type $R_2$, $R_p$ or $R_q$ or
there are proper, nonempty subsets $I$ of $\{0,1,\ldots,p-1\}$ and $J$ of $\{0,1,\ldots,q-1\}$
such that $0\in I$, $I$ has cardinality no more than $(p-1)/2$ and such that a rotation of $h$ is equal to the sorou
\[
\sum_{j\in J^c}\left(\sum_{i\in I}\nu_p^i\right)\nu_q^j+\sum_{j\in J}\left(\sum_{i\in I^c}(-\nu_p^i)\right)\nu_q^j,
\]
where $I^c=\{0,1,\ldots,p-1\}\setminus I$ and $J^c=\{0,1,\ldots,q-1\}\setminus J$.
Thus, $h$ has type
\[
(R_q:\sum_{i\in I}\nu_p^i:|J|R_p).
\]
\end{thm}
\begin{proof}
We may suppose $h$ is not of type $R_2$, $R_p$ of $R_q$.
Thus, (by Lemma~\ref{lem:2p}) its top prime is $q$ and,
after rotating $h$, if necessary, we may write
\[
h=\sum_{j=0}^{q-1}f_j\nu_q^j,
\]
where each $f_j$ is a sum of $2p$-th roots of unity, $\val(f_j)\ne0$ is independent of $j$,
$1\prec f_0$ and $w(f_0)\le w(f_j)$ for all $j$.
Let $a_i$ be the multiplicity of $\nu_p^i$ in $f_0$ and $b_i$ the multiplicity of $-\nu_p^i$ in $f_0$.
Thus,
\[
f_0=\sum_{i=0}^{p-1}a_i\nu_p^i+\sum_{i=0}^{p-1}b_i(-\nu_p^i).
\]
Since $f_0$ has no nonempty vanishing subsorou, we must have $a_ib_i=0$ for all $i$.
Furthermore, letting $I=\{i\in\{0,\ldots,p-1\}\mid a_i>0\}$, we have $0\in I$ and $I$ must be a proper, nonempty subset of $\{0,\ldots,p-1\}$,
for otherwise $f_0$ would have a vanishing subsum of type $R_p$.

Let $j\in\{0,\ldots,1-q\}$.
Since $f_0-f_j$ is a vanishing sum of $2p$-th roots of unity, by Lemma~\ref{lem:2p}, it has type
$R_2^{\oplus k(j)}\oplus R_p^{\oplus \ell(j)}$ for some $k(j),\ell(j)\ge0$.
Though the pair $(k(j),\ell(j))$ is not necessarily unique, we fix a choice, for each $j$.
Moreover, we may (and do) choose $k(0)=w(f_0)$ and $\ell(0)=0$.
We also fix, for each $j$, a decomposition of $f_0-f_j$ into a sum of $k(j)$ minimal vanishing subsorou of type $R_2$
and $\ell(j)$ minimal vanishing subsorou of type $R_p$ and thereby associate, to each term of $f_0$, the corresponding subsorou
of $f_0-f_j$ (either of type $R_2$ or of type $R_p$) in which that term appears.
Of the $a_i$ terms of $f_0$ that are equal to $\nu_p^i$, let $c_i^{(j)}$ be the number of them whose associated subsorou
(as above) of $f_0-f_j$ is of type $R_p$.
Similarly, of the $b_i$ terms of $f_0$ that are equal to $-\nu_p^i$, let $d_i^{(j)}$ be the number of them whose
associated subsorou of $f_0-f_j$ is of type $R_p$.
Note that we have
\begin{equation}\label{eq:ca,db}
0\le c_i^{(j)}\le a_i,\qquad 0\le d_i^{(j)}\le b_i.
\end{equation}

Each of the $\ell(j)$ minimal vanishing subsorou of $f_0-f_j$ of type $R_p$ in the fixed decomposition {of $f_0-f_j$} must be of the form
$1+\nu_p+\nu_p^2+\cdots+\nu_{p}^{p-1}$ or $-1-\nu_p-\nu_p^2-\cdots-\nu_p^{p-1}$.
Let $\ell_+(j)$ and $\ell_-(j)$ be the number of each form, respectively.
Thus, for all $j\in\{0,\ldots,q-1\}$, we have
\begin{equation}\label{eq:ell+-}
\ell_+(j)\ge\max_{0\le i\le p-1}c_i^{(j)},\qquad\ell_-(j)\ge\max_{0\le i\le p-1}d_i^{(j)}.
\end{equation}

We now consider the multiplicity of $\nu_p^i$ in $-f_j$.
There are $b_i-d_i^{(j)}$ terms in $f_0$ equal to $-\nu_p^i$ that appear in one of the $k(j)$ vanishing subsorou of $f_0-f_j$
having type $R_2$ in the fixed decomposition.
These contribute $b_i-d_i^{(j)}$ to the multiplicity of $\nu_p^i$ in $-f_j$.
The other appearances of $\nu_p^i$ in $f_j$ occur in those of the $\ell_+(j)$ subsorou of $f_0-f_j$
of  the form $1+\nu_p+\cdots+\nu_p^{p-1}$
in the fixed decomposition in which the term $\nu_p^i$ does not come from $f_0$.
There are $\ell_+(j)-c_i^{(j)}$ such appearances.
Thus, the multiplicity of $\nu_p^i$ in $-f_j$ is equal to $b_i-d_i^{(j)}+\ell_+(j)-c_i^{(j)}$.
Similarly, the multiplicity of $-\nu_p^i$ in $-f_j$ is $a_i-c_i^{(j)}+\ell_-(j)-d_i^{(j)}$.
Thus, we have
\[
f_j=\sum_{i=0}^{p-1}\big(a_i-c_i^{(j)}+\ell_-(j)-d_i^{(j)}\big)\nu_p^i
+\sum_{i=1}^{p-1}\big(b_i-d_i^{(j)}+\ell_+(j)-c_i^{(j)}\big)(-\nu_p^i).
\]
If $\ell_+(j)=0$, then, by~\eqref{eq:ell+-}, $c_i^{(j)}=0$ for all $i$.
Then, for all $i\in I$, we have
\[
a_i-c_i^{(j)}+\ell_-(j)-d_i^{(j)}\ge a_i>0.
\]
Consequently, we have
\begin{equation}\label{eq:I<fj}
\sum_{i\in I}\nu_p^i\prec f_j.
\end{equation}
If $\ell_+(j)>0$, then since, for all $i\in I^c$ we have $a_i=0$, by~\eqref{eq:ca,db} we have $c_i^{(j)}=0$ and, thus,
\[
b_i-d_i^{(j)}+\ell_+(j)-c_i^{(j)}\ge\ell_+(j)>0.
\]
Consequently, we have
\begin{equation}\label{eq:Ic<fj}
\sum_{i\in I^c}(-\nu_p^i)\prec f_j.
\end{equation}
Let $J=\{j\in\{0,\ldots,q-1\}\mid \ell_+(j)>0\}$.
Since $\ell(0)=0$, we have $0\notin J$.
Using~\eqref{eq:I<fj} and~\eqref{eq:Ic<fj}, we have
\begin{equation}\label{eq:<h}
\sum_{j\in J^c}\left(\sum_{i\in I}\nu_p^i\right)\nu_q^j
+\sum_{j\in J}\left(\sum_{i\in I^c}(-\nu_p^i)\right)\nu_q^j\prec h.
\end{equation}
Since
\[
\val\left(\sum_{i\in I}\nu_p^i\right)=\val\left(\sum_{i\in I^c}(-\nu_p^i)\right),
\]
the sorou on the left-hand-side of~\eqref{eq:<h} vanishes.
Since $h$ is a minimal vanishing sorou, we must have equality in~\eqref{eq:<h}.
The set $J$ cannot be empty, for this would violate that $h$ is minimal vanishing.
Finally, by suitable rotation, the roles of $I$ and $I^c$ can be interchanged, simultaneously
with those of $J$ and $J^c$.
Thus, after rotation if necessary, we may obtain that the cardinality of $I$ is no greater than $(p-1)/2$.
\end{proof}

The following theorem extends the work of Poonen and Rubinstein~\cite{PR98}, whose classification
went up to weight 12.
\begin{thm}\label{thm:mvs}
Types of all
of the minimal vanishing sums of roots of unity of weight no greater than 16, up to rotation, are listed in
Table~\ref{tab:mvs}.
All have height $1$.
Also listed (to help with the derivation) are the top prime
and the weight partition of each type.
Furthermore, the possible parities of orders of terms of the sorou of each type are listed.
All of the indicated possible parities do, in fact, occur.
\end{thm}

\begin{center}
{\footnotesize
\begin{longtable}{c|c|c|c|c|l}
\caption{Minimal vanishing sums of roots of unity and their possible parities of orders.} \\
\label{tab:mvs}
           & Top     & Relative & Weight   &          & Possible\\
Weight & prime & order     & partition & Type & parities  \\ \hline\hline
\endfirsthead
\multicolumn{5}{c}%
{\tablename\ \thetable\ -- \textit{Continued from previous page}} \\
\hline
weight & prime & order &  partition & type & parities \\ \hline
\endhead
\hline \multicolumn{5}{r}{\textit{Continued on next page}} \\
\endfoot
\hline
\endlastfoot
$2$ & $2$ & $2$ & $(1,1)$ & $R_2$ & $(1,1)$ \\ \hline
$3$ & $3$ & $3$ & $(1,1,1)$ & $R_3$ & $(3,0)$ \\ \hline
$5$ & $5$ & $5$ & $(1,1,1,1,1)$ & $R_5$ & $(5,0)$ \\ \hline
$6$ & $5$ & $30$ & $(1,1,1,1,2)$ & $(R_5:R_3)$ & $(4,2)$ \\ \hline\hline
$7$ & $7$ & $7$ & $(1,1,1,1,1,1,1)$ & $R_7$ & $(7,0)$ \\ \cline{2-6}
    & $5$ & $30$ & $(1,1,1,2,2)$ & $(R_5:2R_3)$ & $(4,3)$ \\ \hline\hline
$8$ & $7$ & $42$ & $(1,1,1,1,1,1,2)$ & $(R_7:R_3)$ & $(6,2)$ \\ \cline{2-6}
    & $5$ & $30$ & $(1,1,2,2,2)$ & $(R_5:3R_3)$ & $(6,2)$ \\ \hline\hline
$9$ & $7$ & $42$ & $(1,1,1,1,1,2,2)$ & $(R_7:2R_3)$ & $(5,4)$ \\ \cline{2-6}
    & $5$ & $30$ & $(1,2,2,2,2)$ & $(R_5:4R_3)$ & $(8,1)$  \\ \hline\hline
$10$ & $7$ & $70$ & $(1,1,1,1,1,1,4)$ & $(R_7:R_5)$ & $(6,4)$ \\ \cline{2-6}
    & $7$ & $42$ & $(1,1,1,1,2,2,2)$ & $(R_7:3R_3)$ & $(6,4)$  \\ \hline\hline
$11$ & $11$ & $11$ & $(1,\ldots,1)$ & $R_{11}$ & $(11,0)$ \\ \cline{2-6}
     & $7$ & $210$ & $(1,1,1,1,1,1,5)$ & $(R_7:(R_5:R_3))$ & $(10,1),\,(8,3)$ \\ \cline{2-6}
     & $7$ & $210$ & $(1,1,1,1,1,2,4)$ & $(R_7:R_3,R_5)$ & $(6,5)$ \\ \cline{2-6}
     & $7$ & $42$ & $(1,1,1,2,2,2,2)$ & $(R_7:4R_3)$ & $(8,3)$  \\ \hline\hline
$12$ & $11$ & $66$ & $(1,\ldots,1,2)$ & $(R_{11}:R_3)$ & $(10,2)$ \\ \cline{2-6}
     & $7$ & $210$ & $(1,1,1,1,1,1,6)$ & $(R_7:(R_5:2R_3))$ & $(10,2),\,(9,3)$ \\ \cline{2-6}
     & $7$ & $210$ & $(1,1,1,1,1,2,5)$ & $(R_7:R_3,(R_5:R_3))$ & $(9,3),\,(7,5)$ \\ \cline{2-6}
     & $7$ & $210$ & $(1,1,1,1,2,2,4)$ & $(R_7:2R_3,R_5)$ & $(8,4)$ \\ \cline{2-6}
     & $7$ & $42$ & $(1,1,2,2,2,2,2)$ & $(R_7:5R_3)$ & $(10,2)$  \\ \hline\hline
$13$ & $13$ & $13$ & $(1,\ldots,1)$ & $R_{13}$ & $(13,0)$ \\ \cline{2-6}
     &  $11$ & $66$ & $(1,\ldots,1,2,2)$ & $(R_{11}:2R_3)$ & $(9,4)$ \\ \cline{2-6}
     &  $7$ & $210$ & $(1,1,1,1,1,1,7)$ & $(R_7:(R_5:3R_3))$ & $(12,1),\,(8,5)$ \\ \cline{2-6}
     &  $7$ & $210$ & $(1,1,1,1,1,2,6)$ & $(R_7:R_3,(R_5:2R_3))$ & $(9,4),\,(8,5)$  \\ \cline{2-6}
     &  $7$ & $70$ & $(1,1,1,1,1,4,4)$ & $(R_7:2R_5)$ & $(8,5)$ \\ \cline{2-6}
     &  $7$ & $210$ & $(1,1,1,1,2,2,5)$ & $(R_7:2R_3,(R_5:R_3))$ & $(8,5),\,(7,6)$ \\ \cline{2-6}
     &  $7$ & $210$ & $(1,1,1,2,2,2,4)$ & $(R_7:3R_3,R_5)$ & $(10,3)$ \\ \cline{2-6}
     &  $7$ & $42$ & $(1,2,2,2,2,2,2)$ & $(R_7:6R_3)$ & $(12,1)$ \\  \hline\hline
$14$ & $13$ & $78$ & $(1,\ldots,1,2)$ & $(R_{13}:R_3)$ & $(12,2)$ \\ \cline{2-6}
     &  $11$ & $110$ & $(1,\ldots,1,4)$ & $(R_{11}:R_5)$ & $(10,4)$ \\ \cline{2-6}
     &  $11$ & $66$ & $(1,\ldots,1,2,2,2)$ & $(R_{11}:3R_3)$ & $(8,6)$  \\ \cline{2-6}
     &  $7$ & $210$ & $(1,1,1,1,1,1,8)$ & $(R_7:(R_5:4R_3))$ & $(14,0),\,(7,7)$ \\ \cline{2-6}
     &  $7$ & $210$ & $(1,1,1,1,1,2,7)$ & $(R_7:R_3,(R_5:3R_3))$ & $(11,3),\,(7,7)$ \\ \cline{2-6}
     &  $7$ & $210$ & $(1,1,1,1,1,4,5)$ & $(R_7:R_5,(R_5:R_3))$ & $(9,5),\,(7,7)$ \\ \cline{2-6}
     &  $7$ & $210$ & $(1,1,1,1,2,2,6)$ & $(R_7:2R_3,(R_5:2R_3))$ & $(8,6),\,(7,7)$ \\ \cline{2-6}
     &  $7$ & $210$ & $(1,1,1,1,2,4,4)$ & $(R_7:R_3,2R_5)$ &  $(10,4)$ \\ \cline{2-6}
     &  $7$ & $210$ & $(1,1,1,2,2,2,5)$ & $(R_7:3R_3,(R_5:R_3))$ & $(9,5),\,(7,7)$ \\ \cline{2-6}
     &  $7$ & $210$ & $(1,1,2,2,2,2,4)$ & $(R_7:4R_3,R_5)$ & $(12,2)$  \\ \hline\hline
$15$ & $13$ & $78$ & $(1,\ldots,1,2,2)$ & $(R_{13}:2R_3)$ & $(11,4)$ \\ \cline{2-6}
     &  $11$ & $330$ & $(1,\ldots,1,5)$ & $(R_{11}:(R_5:R_3))$ & $(14,1),\,(12,3)$ \\ \cline{2-6}
     &  $11$ & $330$ & $(1,\ldots,1,2,4)$ & $(R_{11}:R_3,R_5)$ &  $(9,6)$ \\ \cline{2-6}
     &  $11$ & $66$ & $(1,\ldots,1,2,2,2,2)$ & $(R_{11}:4R_3)$ & $(8,7)$ \\ \cline{2-6}
     &  $7$ & $210$ & $(1,1,1,1,1,2,8)$ & $(R_7:R_3,(R_5:4R_3))$ & $(13,2),\,(9,6)$ \\ \cline{2-6}
     &  $7$ & $210$ & $(1,1,1,1,1,4,6)$ & $(R_7:R_5,(R_5:2R_3))$ & $(9,6),\,(8,7)$ \\ \cline{2-6}
     &  $7$ & $210$ & $(1,1,1,1,1,5,5)$ & $(R_7:2(R_5:R_3))$ &  %
$(13,2),\,(11,4),\,(9,6)$ \\      \cline{2-6}
     &  $7$ & $210$ & $(1,1,1,1,2,2,7)$ & $(R_7:2R_3,(R_5:3R_3))$ & $(10,5),\,(9,6)$ \\ \cline{2-6}
     &  $7$ & $210$ & $(1,1,1,1,2,4,5)$ & $(R_7:R_3,R_5,(R_5:R_3))$ & $(9,6),\,(8,7)$ \\ \cline{2-6}
     &  $7$ & $210$ & $(1,1,1,2,2,2,6)$ & $(R_7:3R_3,(R_5:2R_3))$ & $(9,6),\,(8,7)$ \\ \cline{2-6}
     &  $7$ & $210$ & $(1,1,1,2,2,4,4)$ & $(R_7:2R_3,2R_5)$ & $(12,3)$ \\ \cline{2-6}
     &  $7$ & $210$ & $(1,1,2,2,2,2,5)$ & $(R_7:4R_3,(R_5:R_3))$ & $(11,4),\,(9,6)$ \\ \cline{2-6}
     &  $7$ & $210$ & $(1,2,2,2,2,2,4)$ & $(R_7:5R_3,R_5)$ & $(14,1)$ \\ \cline{2-6}
     &  $7$ & $70$ & $(2,2,2,2,2,2,3)$ & $(R_7:1+\nu_5^y:R_5)$ & $(12,3)$ \\*
     &         &          &                  & $y\in\{1,2\}$ & \\ \hline\hline
$16$ & $13$ & $130$ & $(1,\ldots,1,4)$ & $(R_{13}:R_5)$ & $(12,4)$ \\ \cline{2-6}
     &  $13$ & $78$ & $(1,\ldots,1,2,2,2)$ & $(R_{13}:3R_3)$ & $(10,6)$ \\ \cline{2-6}
     &  $11$ & $154$ & $(1,\ldots,1,6)$ & $(R_{11}:R_7)$ &  $(10,6)$ \\ \cline{2-6}
     &  $11$ & $330$ & $(1,\ldots,1,6)$ & $(R_{11}:(R_5:2R_3))$ &  %
$(14,2),\,(13,3)$      \\ \cline{2-6}
     &  $11$ & $330$ & $(1,\ldots,1,2,5)$ & $(R_{11}:R_3,(R_5:R_3))$ & $(13,3),\,(11,5)$ \\ \cline{2-6}
     &  $11$ & $330$ & $(1,\ldots,1,2,2,4)$ & $(R_{11}:2R_3,R_5)$ & $(8,8)$  \\ \cline{2-6}
     &  $11$ & $66$ & $(1,\ldots,1,2,2,2,2,2)$ & $(R_{11}:5R_3)$ & $(10,6)$ \\ \cline{2-6}
     &  $7$ & $210$ & $(1,1,1,1,1,4,7)$ & $(R_7:R_5,(R_5:3R_3))$ & $(11,5),\,(9,7)$ \\ \cline{2-6}
     &  $7$ & $210$ & $(1,1,1,1,1,5,6)$ & $(R_7:(R_5:R_3),(R_5:2R_3))$ & $(13,3),\,(12,4)$ \\*
     &         &            &                &                                            & $(11,5),\,(10,6)$ \\ \cline{2-6}
     &  $7$ & $210$ & $(1,1,1,1,2,2,8)$ & $(R_7:2R_3,(R_5:4R_3))$ & $(12,4),\,(11,5)$ \\ \cline{2-6}
     &  $7$ & $210$ & $(1,1,1,1,2,4,6)$ & $(R_7:R_3,R_5,(R_5:2R_3))$ & $(9,7),\,(8,8)$ \\ \cline{2-6}
     &  $7$ & $210$ & $(1,1,1,1,2,5,5)$ & $(R_7:R_3,2(R_5:R_3))$ & $(12,4),\,(10,6),\,(8,8)$ \\ \cline{2-6} %
     &  $7$ & $210$ & $(1,1,1,1,4,4,4)$ & $(R_7:3R_5)$ & $(12,4)$ \\ \cline{2-6} %
     &  $7$ & $210$ & $(1,1,1,2,2,2,7)$ & $(R_7:3R_3,(R_5:3R_3))$ & $(11,5),\,(9,7)$ \\ \cline{2-6}
     &  $7$ & $210$ & $(1,1,1,2,2,4,5)$ & $(R_7:2R_3,R_5,(R_5:R_3))$ & $(11,5),\,(9,7)$ \\ \cline{2-6}
     &  $7$ & $210$ & $(1,1,2,2,2,2,6)$ & $(R_7:4R_3,(R_5:2R_3))$ & $(11,5),\,(10,6)$ \\ \cline{2-6}
     &  $7$ & $210$ & $(1,1,2,2,2,4,4)$ & $(R_7:3R_3,2R_5)$ & $(14,2)$ \\ \cline{2-6}
     &  $7$ & $210$ & $(1,2,2,2,2,2,5)$ & $(R_7:5R_3,(R_5:R_3))$ & $(13,3),\,(11,5)$ \\ \cline{2-6}
     &  $7$ & $105$ & $(2,2,2,2,2,2,4)$ & $(R_7:1+\nu_3:(R_5:R_3))$ & $(16,0)$ \\ \cline{2-6}
     &  $7$ & $210$ & $(2,2,2,2,2,2,4)$ & $(R_7:1+\nu_5^y:(R_5:R_3))$ & $(14,2)$ \\*
     &         &           &                    & $y\in\{1,2\}$ &  \\ \cline{2-6}
     &  $7$ & $210$ & $(2,2,2,2,2,2,4)$ & $(R_7:1-\nu_3\nu_5^y:(R_5:R_3))$ & $(9,7)$ \\*
     &         &           &                    & $y\in\{1,2,3,4\}$ &  \\ \cline{2-6}
     &  $7$ & $210$ & $(2,2,2,2,2,3,3)$ & $(R_7:1+\nu_5^y:R_2\oplus R_3,\,R_5)$ & $(11,5)$ \\*
     &         &            &                    & $y\in\{1,2\}$ &  \\ \cline{2-6}
     &  $7$ & $70$ & $(2,2,2,2,2,3,3)$ & $(R_7:1+\nu_5^y:2R_5)$ & $(10,6)$ \\*
     &         &          &                    & $y\in\{1,2\}$ &  \\ \hline
\end{longtable} }
\end{center}

{Regarding parities, we have} for example, {that} a minimal vanishing sorou of type $(R_5:R_3)$ can have either 4 terms of even order and 2 of odd order,
or 4 of odd order and 2 of even order, while one of type $(R_7:(R_5,R_3))$ can have either 10 even and 1 odd, 8 even and 3 odd,
3 odd and 8 even or 1 odd and 10 even.

\begin{proof}[Proof of Theorem~\ref{thm:mvs}]
First of all, note that, by Theorem~\ref{thm:2pq}, if the top prime of a minimal vanishing sorou is $\le 5$, then
the weight partition must contain $1$.

For those types whose weight partition contains $1$, the cataloging is essentially a continuation of
the method of Poonen and Rubinstein~\cite{PR98}.
We will now describe the general technique in this case, while we work through two examples.

Suppose $h$ is a minimal vanishing sorou with weight $n$.
For example, say $n=13$.
Let $p$ be the largest prime divisor of the relative order of $h$.
For example, say $p=7$.
Write $h$ as in equation~\eqref{eq:hnuf} of Proposition~\ref{prop:hnuf}.
Then $(w(f_0),\ldots,w(f_{p-1}))$ is a partition of $n$ into $p$ parts,
the nondecreasing re-ordering of which will be the weight partition.
For this portion of the proof, we are supposing that $1$ appears in this partition.
After rotating $h$, if necessary, we may without loss of generality assume $f_0=1$.
For example, suppose the weight partition is $(1,1,1,1,1,4,4)$.
By Proposition~\ref{prop:hnuf},
whenever $w(f_j)=1$ but $j>0$, we need $f_0-f_j$ to be vanishing of weight $2$, so we must have $f_j=1$.
If $w(f_j)>1$, then we need $f_0-f_j$ to be vanishing.
Moreover, since $f_0=1$ and $f_j$ must have no vanishing nonempty subsorou, $f_0-f_j$ must be a minimal vanishing sorou
of weight $w(f_j)+1$ and must have $1$ as a subsorou and must have relative order that is a product of primes that are each
strictly less than $p$.
Thus, the construction of this portion of the table (when the weight partition contains $1$) proceeds recursively.
In the case of $p=7$ and $w(f_j)=4$, we must have $f_0-f_j$ of type $R_5$.
This forces $f_0-f_j$ to have type $R_5$ and $f_j$ to be $-\nu_5-\nu_5^2-\nu_5^3-\nu_5^4$.
This implies that the order parities are $(8,5)$ --- in this case, five odd order terms
of the form $\nu_7^j$ for the $j$ such that $f_j=1$,
and a total of eight even order terms {from two sums} of the form $\nu_7^j(-\nu_5-\nu_5^2-\nu_5^3-\nu_5^4)$ for the $j$ such that $w(f_j)=4$.
After rotation, of course, the odds and the evens can be interchanged.
This is the reason we don't specify how many of each on the ``Parities'' column of Table~\ref{tab:mvs}.
This shows that all minimal vanishing sorou with prime $7$ and weight partition $(1,1,1,1,1,4,4)$ have {parities $(8,5)$.}
Conversely, using Proposition~\ref{prop:hnuf}, we see that minimal vanishing sorou with this prime and weight partition do occur.

We now consider the case of weight partition $(1,1,1,1,2,2,5)$.
Arguing as before, if $h$ is a vanishing sorou with this weight partition, after rotating we have $f_0=1$ and, for those
$j$ with $w(f_j)=1$ we have $f_j=1$.
For those $j$ with $w(f_j)=2$, since $f_0-f_j$ is a minimal vanishing sorou, we need $f_j=-\nu_3-\nu_3^2$.
For $j$ with $w(f_j)=5$, $f_0-f_j$ must have type $(R_5,R_3)$, meaning it must be equal to a rotation of
\[
1+\nu_5+\nu_5^2+\nu_5^3+\nu_5^4(-\nu_3-\nu_3^2).
\]
Thus, $f_j$ must be one of the sorou
\begin{gather*}
-\nu_5-\nu_5^2-\nu_5^3+\nu_5^4(\nu_3+\nu_3^2), \\
-\nu_5-\nu_5^2+\nu_5^3(\nu_3+\nu_3^2)-\nu_5^4, \\
-\nu_5+\nu_5^2(\nu_3+\nu_3^2)-\nu_5^3-\nu_5^4, \\
\nu_5(\nu_3+\nu_3^2)-\nu_5^2-\nu_5^3-\nu_5^4, \\
-\nu_3+\nu_3^2(\nu_5+\nu_5^2+\nu_5^3+\nu_5^4) \\
\nu_3(\nu_5+\nu_5^2+\nu_5^3+\nu_5^4) -\nu_3^2
\end{gather*}
Counting the terms with odd and, respectively, even orders, we get (unordered)
parity distributions $(8,5)$ and $(7,6)$ {for $h$}, depending on which of the above six is chosen.
Conversely, we see that each choice leads to a minimal vanishing sorou with the given weight partition, {so} that both
parity distributions occur.

We now analyze the situations where the smallest subsidiary weight of a minimal vanishing sorou $h$ is strictly greater than $1$.
By Theorem~\ref{thm:2pq}, this is possible only for top prime at least $7$.
We cannot have weight partition of the form $(2,2,\ldots,2)$, for this would require all subsidiary types to be $R_2\oplus R_2$,
which would imply, writing $h$ as in~\eqref{eq:hnuf}, that $f_0=f_1=\cdots=f_{p-1}$ and then $h$ would have a proper, nonempty vanishing subsorou
of type $R_p$.
More generally, if the smallest subsidiary weight is equal to $2$, then
{assuming (without loss of generality) $1\prec f_0$ and} letting
$(R_p:1+b:T_1,\ldots,T_n)$ be the type of $h$,
at least one of the subsidiary types $T_1,\ldots,T_n$ must be minimal vanishing.
Indeed if not, then we must be able to write each $f_0-f_j=g_{j,1}+g_{j_2}$ for $g_{j,i}$ vanishing of type $T_{j,i}$, with
$T=T_{j,1}\oplus T_{j,2}$.
Then, after interchanging $g_{j,1}$ and $g_{j,2}$, if necessary, must have $1\prec g_{j,1}$ and $b\prec g_{g,2}$,
for otherwise we would get a proper, vanishing subsorou $\nu_p^j g_{j,i}$ of $h$.
But then, separating the bits, this implies that $h$ has a proper, vanishing subsorou of type
\begin{equation}\label{eq:RTj1}
(R_p:T_{1,1},T_{2,1},\ldots,T_{n,1}).
\end{equation}

Thus, the {lowest weight} possibility is with top prime $7$, weight partition $(2,2,2,2,2,2,3)$ and total weight $15$.
Suppose $h$ is a minimal vanishing sorou with this data.
Then
\begin{equation}\label{eq:h7}
h=\sum_{j=0}^6\nu_7^jf_j
\end{equation}
and after rotation, we may assume $1\prec f_0$, $w(f_0)=w(f_1)=\cdots=w(f_5)=2$ and $w(f_6)=3$.
{Hence,} $f_0=1+b$ with $b$ a root of unity whose order divides $30$.
For $1\le j\le 5$, $f_0-f_j$ must have type $R_2\oplus R_2$ and, therefore, we have $f_j=f_0$.
The type of $f_0-f_6$ is either $R_2\oplus R_3$ or $R_5$.
{Since,} by the earlier remark, {at least one of the subsidiary types must be minimal vanishing, the type of $f_0-f_6$ must be} $R_5$.
After rotation, we must have $b\in\{\nu_5,\nu_5^2\}$.
These yield the two minimal vanishing sorou
\begin{gather*}
\left(\sum_{j=0}^6\nu_7^j(1+\nu_5)\right)+\nu_7^6(-\nu_5^2-\nu_5^3-\nu_5^4) \\
\left(\sum_{j=0}^6\nu_7^j(1+\nu_5^2)\right)+\nu_7^6(-\nu_5-\nu_5^3-\nu_5^4), 
\end{gather*}
having height $1$ and parities $(12,3)$.
These are of order $2\cdot5\cdot7$ and were, in fact, treated in Theorem~\ref{thm:2pq}.

We now turn to the case of top prime $p=7$ and weight partition $(2,2,2,2,2,2,4)$, for total weight $16$.
Arguing as before, after rotation we must have $h$ as in~\eqref{eq:h7} with $f_0=f_1=\cdots=f_5=1+b$ and with $f_0-f_6$
vanishing of weight $6$.
By the observation made above~\eqref{eq:RTj1}, $f_0-f_6$ must be minimal vanishing of weight $6$, so
must be of type $(R_5:R_3)$.
Since $1\prec f_0$, after rotation we must have $f_0=1+\nu_5^y$,
some $y\in\{1,2\}$, or $f_0=1+\nu_3$
or $f_0=1-\nu_3\nu_5^y$, some $y\in\{1,2,3,4\}$.
In the first case, when $f_0=1+\nu_5^y$,
there are three different choices for $f_6$, but they all yield minimal vanishing sorou of height $1$ and parities $(14,2)$.
For example, when $f_0=1+\nu_5$, then we must have either
\[
{h=\sum_{j=0}^5\nu_7^j(1+\nu_5)+\nu_7^6(-\nu_5^2-\nu_5^3+\nu_5^4(\nu_3+\nu_3^2))}
\]
{or one of two other similar possibilities.}
When $f_0=1+\nu_3$, then we must have $f_6=\nu_3^2(\nu_5+\nu_5^2+\nu_5^3+\nu_5^4)$, and we get
\[
h=\sum_{j=0}^5\nu_7^j(1+\nu_3)+\nu_7^6\nu_3^2(\nu_5+\nu_5^2+\nu_5^3+\nu_5^4)
\]
with height $1$ and parities $(16,0)$.
The third case yields {$h$ with height $1$ and parities $(9,7)$.}
For example if $f_0=1-\nu_3\nu_5$, the $f_6$ must be one of
\[
\nu_3^2\nu_5-\nu_5^2-\nu_5^3-\nu_5^4\quad\text{or}\quad
-\nu_3^2+\nu_3\nu_5^2+\nu_3\nu_5^3+\nu_3\nu_5^4,
\]
{yielding
\begin{align*}
h&=\sum_{j=0}^5\nu_7^j(1-\nu_3\nu_5)+\nu_7^6(\nu_3^2\nu_5-\nu_5^2-\nu_5^3-\nu_5^4), \\
h&=\sum_{j=0}^5\nu_7^j(1-\nu_3\nu_5)+\nu_7^6(-\nu_3^2+\nu_3\nu_5^2+\nu_3\nu_5^3+\nu_3\nu_5^4),
\end{align*}
respectively.}

The next possibility is with top prime $7$, weight partition $(2,2,2,2,2,3,3)$ and total weight $16$.
Suppose $h$ is a minimal vanishing sorou with this data.
Then $h$ is as in~\eqref{eq:h7}
and after rotation, we may assume $1\prec f_0$ and for a set $J=\{j(1),j(2)\}$ equal to one of $\{5,6\}$,
$\{4,6\}$ or $\{3,6\}$, we have
\[
w(f_j)=\begin{cases}3,&j\in J,\\2,&j\notin J.\end{cases}
\]
We have $f_0=1+b$ with $b$ a root of unity not equal to $-1$ and having order dividing $30$.
As argued above, if $j\notin J$, then $f_j=f_0$.
For $j\in J$,
the type of $f_0-f_j$ is either $R_2\oplus R_3$ or $R_5$.
By the observation made above~\eqref{eq:RTj1}, they cannot for both $j\in J$ be of type $R_2\oplus R_3$.

If both types $R_2\oplus R_3$ and $R_5$ occur, then without loss of generality,
the type of $f_0-f_{j(1)}$ is $R_2\oplus R_3$ and the type of $f_0-f_{j(2)}$ is $R_5$.
{After} rotation,
we must have $f_0=1+\nu_5^y$ for some $y\in\{1,2\}$.
$f_{j(2)}=-\nu_5^{3-y}-\nu_5^3-\nu_5^4$
and $f_{j(1)}$ equal to either $1+\nu_5^y(-\nu_3-\nu_3^2)$ or $\nu_5^y-\nu_3-\nu_3^2$.
Thus, we get the {two} minimal vanishing sorou
\begin{align*}
{h}&=\left(\sum_{j\in J^c}\nu_7^j(1+\nu_5^y)\right)+\nu_7^{j(1)}\big(1+\nu_5^y(-\nu_3-\nu_3^2)\big)+\nu_7^{j(2)}(-\nu_5^{3-y}-\nu_5^3-\nu_5^4), \\
{h}&=\left(\sum_{j\in J^c}\nu_7^j(1+\nu_5^y)\right)+\nu_7^{j(1)}(\nu_5^y-\nu_3-\nu_3^2)+\nu_7^{j(2)}(-\nu_5^{3-y}-\nu_5^3-\nu_5^4),
\end{align*}
{each} of type
\[
(R_7:1+\nu_5^y:R_2\oplus R_3,\,R_5),
\]
height $1$ and with parities $(11,5)$.

The case when the type of $f_0-f_j$ is $R_5$ for both $j\in J$ is treated similarly.
We get minimal vanishing sorou of type $(R_7:1+\nu_5^y:2R_5)$, all with height $1$, order $70=2\cdot5\cdot7$
and parities $(10,6)$.
{These cases were} also treated in Theorem~\ref{thm:2pq}.

Regarding heights, if a minimal vanishing sorou $h$ has type $(R_p:f_0:T_1,\ldots,T_n)$, then, writing $h$ as in~\eqref{eq:hnuf},
we see that the height of $h$ is the maximum of the heights of the $f_j$ for $0\le j\le p-1$.
Thus, the height of $h$ is no greater than the maximum of the heights of $f_0$ and of the partitions $f_0-f_j$ (when $f_j\ne f_0$)
whose types are the subsidiary types $T_1,\ldots,T_n$.
In particular, $h$ can have height strictly greater than $1$ only if either $f_0$ has height $>1$ or some vanishing sorou
of type $T_i$ has height $>1$, for some $i$.
Examining Table~\ref{tab:mvs}, we see that for weight $\le15$, all the $f_0$ have height $1$ and all the subsidiary types are of minimal vanishing
sorou.
Thus, by induction, all minimal vanishing sorou of weight $\le15$ have height $1$.

The first example of a minimal vanishing sorou with non-minimal vanishing subsidiary types is 
the case of those of type $(R_7:1+\nu_5^y:R_2\oplus R_3,\,R_5)$, of weight $16$.
However, the analysis of these conducted above shows that they all have height $1$.
\end{proof}

\section{Towards a classification of vanishing sorou of weight $>16$}\label{sec:comp1}

To classify sorou of a large weight, computer assistance becomes indispensable.
We designed an algorithm that
is essentially an exhaustive search, 
which we are able to do because results from previous sections
imply that the search space is finite.

Some preliminary  code  in python 
is available from the {website \url{https://github.com/lchristie/Sums-of-Roots-of-Unity}.}
The main part of the algorithm is contained in {\tt TypeGen.py} and described in {Algorithm \ref{algo:TypeGen}.}

\begin{remark}[Computational nature] The data and Conjecture \ref{conj:ht2} below are based on a preliminary Python implementation of Algorithm \ref{algo:TypeGen}. 
{The code used for the computations in this section was written manually (without large language models or other AI-assisted code generation tools).}
We have not formally verified the correctness or completeness of this implementation, so these results should be regarded as experimental evidence.
{The computational results also rely on floating-point vanishing tests, so numerical error is another potential source of uncertainty.}
\end{remark}

{
\begin{remark}[Theoretical scope of the algorithms]
For each fixed integer $N$, Algorithms \ref{algo:1} and \ref{algo:2}  below describe an exhaustive \emph{finite} search for minimal vanishing sorou of weight $\le N$ (and hence of all possible heights, since the height of a sorou of weight $N$ is automatically $\le N$). Indeed, the weight bound forces a finite set of candidate top primes $p\le N$ and a finite set of possible $f_0$'s, and the recursive generation of sorou terminates because each recursive call strictly decreases the weight. Consequently, assuming a correct implementation (and an exact vanishing test), these algorithms yield a complete classification for each fixed weight bound $N$, and hence an enumeration of all minimal vanishing sorou by increasing $N$. 
\end{remark}
}

\begin{algo}\label{algo:TypeGen}\rm
Algorithm 1 below uses the existing classification of 
minimal vanishing sorou up to weight $\leq k$ to find all minimal vanishing sorou of weight $k+1$.
We use the classification obtained in Section \ref{sec:hand} to start if off.
{(We could have started this algorithm using {the} list~\cite{PR98} of minimal vanishing sorou up to weight $12$,
due to 
Poonen and Rubinstein.
However the derivation presented in Section~\ref{sec:hand} of minimal vanishing sorou up to weight $16$ was important in order to develop
understanding of  and to present the expanded notion of type.)}\looseness=-1

We now describe how the algorithm finds all minimal vanishing types of weight $k+1$.
First note that 
the only possible top primes 
for a vanishing sorou of weight $k+1$ are $\leq k+1$.
For each of these primes $p$ we build the possible partitions
of $k+1$ into $p$ summands. 
For each of these partitions we generate all possible $f_0$'s
(cf.~Definition \ref{defi:type}). 
Note that $f_0$ will be a sorou of relative order 
dividing the product of primes $<p$; this observation together 
with the bound $k+1$ on the weight makes this a finite list.

Then for each of these $f_0$ we first prepare 
all possible lists of subtypes with weights given by the
partition (this uses the already constructed type list of weight $<k+1$). 
This is done carefully to eliminate redundancy and repetition.
For each list we then build a new type of weight $k+1$, generate a single sorou 
of this type and use this to check whether it is minimal vanishing with 
Proposition \ref{prop:hnuf}.
{To test whether a sorou 
$h$ vanishes, we evaluate 
$\val(h)$ numerically (floating-point complex arithmetic) and declare it vanishing when 
$|\val(h)|$ is below a fixed tolerance. Exact algebraic-number arithmetic in SymPy is in principle possible but was too slow for the exhaustive search.}

It is now straightforward to verify
that this algorithm generates all minimal vanishing types of weight $k+1$.

\begin{algorithm}
\caption{Generating all minimal vanishing types}\label{algo:1}
\begin{algorithmic}[1]
\Procedure{GenNextTypes}{$previousTypes$}
	\State $Output \leftarrow \varnothing$, $previousWeight \leftarrow \max \{ w(T) : T \in previousTypes\}$
	\State $w_0 \leftarrow previousWeight + 1$ \Comment{Current weight to generate}
	\State $primesToCheck \leftarrow \{ p \in P: p \leq w_0 \}$
	
	\For{$p$ in $primesToCheck$}
		\State $allParitions \leftarrow \{ (x_0, x_1, \dots, x_{p-1}) \in \mathbb{Z}_+^p : \sum x_i = w_0, 0 < x_i \leq x_j \text{ if } i \leq j \}$
		
		\For{$x$ in $allParitions$}
			\If{ $x_0 = \cdots = x_{p-1} = 1 $} \Return $R_p$
			\Else
			\State $F \leftarrow \{ f \in S : w(f) = x_0, relOrder(f) | \prod_{q \in P}^{p-1} q \}$
			
			\For{$f_0$ in $F$}
				\State $X \leftarrow \mathbf{T}_{x_0 + x_0} \times \cdots \times \mathbf{T}_{x_{p-1} + x_0 }$ 	\Comment{Possible subsidiary types}
				\State Filter $X$, based on the conditions below.
				\For{$\mathbf{x}$ in $X$}
					\State Generate one sorou $h$ of type $T = (R_p : f_0 : x_0, \dots, x_{p-1})$. 
					\If{$h$ is minimal vanishing} \Comment{By Proposition 2.3}
					\State $Output \leftarrow Output \cup \{T\}$ \EndIf
				\EndFor
			\EndFor
			\EndIf
		\EndFor
	\EndFor
	\State \Return $Output$ 
\EndProcedure
\end{algorithmic}
\end{algorithm}
\end{algo}

\subsection{Heights and \esorou}\label{sec:comp2}

{Our aforementioned interest in minimal vanishing sorou was actually in
{\em equisigned} ones,}
namely, those whose parities $(n_{\text{odd}},n_{\text{even}})$ satisfy $n_{\text{odd}}=n_{\text{even}}$.
Given the recursive algorithm to describe types, {in order to catalog the equisigned ones, it was} necessary to keep track of the
possible parities for all types of sorou.
To find all parities of a given type we can again use computers. Our preliminary  python code is available from \url{https://github.com/lchristie/Sums-of-Roots-of-Unity}.
The algorithm is contained in {\tt sorouGenerator.py} and described in {Algorithm \ref{algo:parities}.}

\begin{algo}\label{algo:parities}\rm
To find all the parities we need to be able to find all sorou (up to rotation) of a given minimal
vanishing type. We do this recursively per Algorithm 2 below, using the procedure \textit{GenSorou}. %
For a given type $(R_p:f_0:T_1,\ldots,T_n)$ we find all sorou of types $T_i$.
We take combinations of one sorou from each $T_i$, and for each of these combinations we 
subtract all suitable rotations of each sorou from the base sorou $f_{base} = \sum_{i=0}^{p-1} \nu_p^i f_0$ (i.e., rotations by $\nu_q^j$ where $q$ is the relative order of the sorou, $f_0$, and $\nu_p$ concatenated as a formal sorou and $j \in \{0, 1,\dots, q-1\}$, which cancel all terms of $f_0\nu_p^i$).

However,  $T_1,\ldots,T_n$ are not necessarily minimal vanishing
(cf.~Remark \ref{rem:type}\,(f)). We thus also need to generate
every sorou of a given \emph{non-minimal} vanishing type $S_1\oplus\cdots\oplus S_k$.
But it must contain {at least one term from} $f_0$ {in each minimal vanishing component},
which makes the search space finite. We split $f_0$ into $k$ subsorou
and match sorou from each $S_i$ to these subsorou. These sums of the sorou 
are then the possible sorou of type $S_1\oplus\cdots\oplus S_k$ containing $f_0$. This is done via the procedure \textit{GenNonMinSorou}. 

This means that the procedures \textit{GenSorou} and \textit{GenNonMinSorou} recursively call each other. Note however that this is guaranteed to terminate as each call forces a strict reduction in weight. Moreover, types with non-minimal subtypes only occur at weight 16 and higher, so for all of the computations done so far \textit{GenNonMinSorou} never calls itself. 

Once all the sorou of a given type have been constructed, it is easy to check their parities
and heights.
\begin{algorithm}
\caption{Generating All Possible Sorou of a Type}\label{algo:2}
\begin{algorithmic}[1]
	\Procedure{GenSorou}{$T = (R_p : f_0 : T_1, \dots, T_n)$}
	\State $S_0 \leftarrow \{ 0 \}, Output \leftarrow \varnothing$
	\For{$i$ in $\{ 1, \dots, n \}$}	
		\State $S_i \leftarrow GenNonMinSorou(T_i, f_0)$ \Comment{As $T_i$ may not be minimal}
	\EndFor
	\State $\Pi \leftarrow$ All Permutations of $(1, \dots, n, 0, \dots, 0)$ \Comment{With $p - n > 0$ zeros}
	\For{$\pi$ in $\Pi$}
		\State $\Sigma_\pi \leftarrow S_{\pi(1)} \times \cdots \times S_{\pi(p)}$
		\State $f_{base} \leftarrow \sum_{i = 0}^{p-1} f_0 \nu_p^i$
		\For{$(h_1, \dots, h_p)$ in $\Sigma_\pi$}
			\For{$i$ in $\{ 1, \dots p \}$}
			\State $R_i \leftarrow \{ \nu_a^b :  f_0 \preceq \nu_a^b h_i \}$
			\EndFor
			\For{$(\nu_{a_1}^{b_1}, \dots, \nu_{a_p}^{b_p})$ in $R_1 \times \cdots \times R_p$}
				\State $g \leftarrow f_{base} - \sum_{i = 1}^p \nu_p^i \nu_{a_i}^{b_i} h_i$ \Comment{With subtracted terms removed}
				\State $Output \leftarrow Output \cup \{ g \}$.
			\EndFor
		\EndFor
	\EndFor
	\State \Return $Output$
	\EndProcedure
\Procedure{GenNonMinSorou}{$T = T_1 \oplus \cdots \oplus T_m$, $f$}
	\If{$m=1$} \Return $GenSorou(T)$ 
	\Else 
	\State $Output \leftarrow \varnothing$
	\For{$i$ in $\{ 1, \dots, m \}$} 
		\State $S_i \leftarrow GenSorou(T_i)$
	\EndFor
	\State $\mathcal{P} \leftarrow \{ ( h_1, \dots, h_m) : f_0 = \sum_{i= 1}^m h_i, w(h_i) > 0 \}$ \Comment{``Partitions'' of $f_0$}
	\State $\Pi \leftarrow$ All permutations of $\{ 1 , \dots, m \}$
	\State $X \leftarrow \{ (h, \pi, s) \in \mathcal{P} \times  \Pi \times ( S_1 \times \cdots \times S_m) : \forall i \exists \nu_i {}_\ni  h_{\pi(i)} \preceq \nu_i s_i \}$
	\For{$x \in X$}
		\State $R_x \leftarrow \{ (\nu_{a_1}^{b_1}, \dots, \nu_{a_m}^{b_m}) : h_{\pi(i)} \preceq \nu_{a_i}^{b_i} s_i \}$
		\For{ $(\nu_{a_1}^{b_1}, \dots, \nu_{a_m}^{b_m})$ in $R_x$}
			\State $g \leftarrow \sum_{i = 1}^m \nu_{a_i}^{b_i} s_i$ \Comment{So that $f_0 \preceq g$}
			\State $Output \leftarrow Output \cup \{ g \}$
		\EndFor
	\EndFor
	\State \Return $Output$
	\EndIf
	
\EndProcedure
\end{algorithmic}
\end{algorithm}
\end{algo}

\subsection{Computational results and conjectures.}

Appendix \ref{app:all} contains a conjectured complete list of all minimal vanishing types of length $\leq 21$ produced
by 
our implementation of 
this algorithm, as well as their heights and parities. In particular, we conjecture:

\begin{conj}\label{conj:ht2}\mbox{}\par
\begin{enumerate}[label={\rm(\alph*)}]
\item
All minimal vanishing sorou having weight less than 21
are of height $1$. 
\item
{There are minimal vanishing sorou of weight 21 and height strictly greater than $1$; they all have height $2$.}
\end{enumerate}
\end{conj}

\begin{remark}
Minimal vanishing sorou of arbitrarily large
height have been constructed by 
Steinberger \cite{Ste08},
and their existence
can be deduced from Schur's 1931 result
(see \cite{Leh36})
that cyclotomic polynomials $\Phi_n$
have unbounded coefficients.
We point out that 
\begin{multline*}
\Phi_{105}(x)=x^{48}+x^{47}+x^{46}-x^{43}-x^{42}-2
   x^{41}-x^{40}-x^{39}+x^{36}+x^{35}+x^{34}+x^{33}+x^{32}+x^{31}-x^{28}\\ -x^{26} -x^{24}-x^
   {22}-x^{20}+x^{17}+x^{16}+x^{15}+x^{14}+x^{13}+x^{12}-x^9-x^8-2 x^7-x^6-x^5+x^2+x+1
 \end{multline*}
 is the first cyclotomic polynomial
with a coefficient different from $\pm1$.
{It corresponds to a sorou of height $2$ and weight $35$.}
\end{remark}

\newpage
\appendix

\section{Conjectured list of all minimal vanishing types of length $\leq 21$}\label{app:all}

\begin{center}
{\footnotesize
\begin{longtable}{c|c|c|l}
\caption{Conjectured list of minimal vanishing sums of roots of unity and their possible heights and parities of orders.} \\
\label{tab:mvs21}
       &      & Possible & Possible\\
Weight & Type & heights  & parities  \\ \hline\hline
\endfirsthead
\multicolumn{3}{c}%
{\tablename\ \thetable\ -- \textit{Continued from previous page}} \\
\hline
Weight & Type & Heights & Parities \\ \hline
\endhead
\hline \multicolumn{3}{r}{\textit{Continued on next page}} \\
\endfoot
\hline
\endlastfoot

$ 2 $ & $  R_{2}  $ & $  1  $ & $  (1,1)  $ \\ \hline
$ 3 $ & $  R_{3}  $ & $  1  $ & $  (3,0)  $ \\ \hline
$ 5 $ & $  R_{5}  $ & $  1  $ & $  (5,0)  $ \\ \hline
$ 6 $ & $  ( R_{5} : R_{3} )  $ & $  1  $ & $  (4,2)  $ \\ \hline
$ 7 $ & $  ( R_{5} : 2R_{3} ) $ & $  1  $ & $  (4,3)  $ \\ \hline
$ 7 $ & $  R_{7}  $ & $  1  $ & $  (7,0)  $ \\ \hline
$ 8 $ & $  ( R_{5} : 3R_{3} ) $ & $  1  $ & $  (6,2)  $ \\ \hline
$ 8 $ & $  ( R_{7} : R_{3} )  $ & $  1  $ & $  (6,2)  $ \\
\hline

\end{longtable} }
\end{center}


\begin{thebibliography}{9}

\bib{AS12}{article}{
    AUTHOR = {Aliev, I.},
    author={Smyth, C.},
     TITLE = {Solving algebraic equations in roots of unity},
   JOURNAL = {Forum Math.},
    VOLUME = {24},
      YEAR = {2012},
    NUMBER = {3},
     PAGES = {641--665},
      ISSN = {0933-7741},
       URL = {https://doi.org/10.1515/form.2011.087},
}
	
\bib{CW11}{article}{
    AUTHOR = {Chien, T.-Y.},
    author={Waldron, S.},
     TITLE = {A classification of the harmonic frames up to unitary
              equivalence},
   JOURNAL = {Appl. Comput. Harmon. Anal.},
    VOLUME = {30},
      YEAR = {2011},
    NUMBER = {3},
     PAGES = {307--318},
      ISSN = {1063-5203},
       URL = {https://doi.org/10.1016/j.acha.2010.09.001},
}

\bib{CJ76}{article}{
    AUTHOR = {Conway, J. H.},
    author={Jones, A. J.},
     TITLE = {Trigonometric {D}iophantine equations ({O}n vanishing sums of
              roots of unity)},
   JOURNAL = {Acta Arith.},
    VOLUME = {30},
      YEAR = {1976},
    NUMBER = {3},
     PAGES = {229--240},
      ISSN = {0065-1036},
         URL = {https://doi.org/10.4064/aa-30-3-229-240},
}

\bib{DZ00}{article}{
    AUTHOR = {Dvornicich, R.},
    author={Zannier, U.},
     TITLE = {On sums of roots of unity},
   JOURNAL = {Monatsh. Math.},
    VOLUME = {129},
      YEAR = {2000},
    NUMBER = {2},
     PAGES = {97--108},
      ISSN = {0026-9255},
       URL = {https://doi.org/10.1007/s006050050009},
}


\bib{FG14}{article}{
author={Fonseca, R. M.},
author={Grimus, W.},
title={Classification of lepton mixing matrices from finite residual symmetries},
journal={Journal of High Energy Physics},
year={2014},
volume={2014},
number={9},
pages={33 pages},
issn={1029-8479},
url={https://doi.org/10.1007/JHEP09(2014)033}
}



\bib{JMS14}{article}{
    AUTHOR = {Jones, V. F. R.},
    author={Morrison, S.},
    author={Snyder, N.},
     TITLE = {The classification of subfactors of index at most 5},
   JOURNAL = {Bull. Amer. Math. Soc. (N.S.)},
    VOLUME = {51},
      YEAR = {2014},
    NUMBER = {2},
     PAGES = {277--327},
      ISSN = {0273-0979},
       URL = {https://doi.org/10.1090/S0273-0979-2013-01442-3},
}
  



\bib{LL00}{article}{
    AUTHOR = {Lam, T. Y.},
    author={Leung, K. H.},
     TITLE = {On vanishing sums of roots of unity},
   JOURNAL = {J. Algebra},
    VOLUME = {224},
      YEAR = {2000},
    NUMBER = {1},
     PAGES = {91--109},
      ISSN = {0021-8693},
       URL = {https://doi.org/10.1006/jabr.1999.8089},
}



\bib{Leh36}{article}{
    AUTHOR = {Lehmer, E.},
     TITLE = {On the magnitude of the coefficients of the cyclotomic
              polynomial},
   JOURNAL = {Bull. Amer. Math. Soc.},
    VOLUME = {42},
      YEAR = {1936},
    NUMBER = {6},
     PAGES = {389--392},
      ISSN = {0002-9904},
}
		

\bib{Len79}{article}{
    AUTHOR = {{Lenstra, Jr.}, H. W.},
     TITLE = {Vanishing sums of roots of unity},
 BOOKTITLE = {Proceedings, {B}icentennial {C}ongress {W}iskundig
              {G}enootschap ({V}rije {U}niv., {A}msterdam, 1978), {P}art
              {II}},
    SERIES = {Math. Centre Tracts},
    VOLUME = {101},
     PAGES = {249--268},
 PUBLISHER = {Math. Centrum, Amsterdam},
      YEAR = {1979},
}

\bib{Man65}{article}{
   author={Mann, Henry B.},
   title={On linear relations between roots of unity},
   journal={Mathematika},
   volume={12},
   date={1965},
   pages={107--117},
}

\bib{PR98}{article}{
   author={Poonen, Bjorn},
   author={Rubinstein, Michael},
   title={The number of intersection points made by the diagonals of a
   regular polygon},
   journal={SIAM J. Discrete Math.},
   volume={11},
   date={1998},
   pages={135--156},
}

\bib{Sch64}{article}{
    AUTHOR = {Schoenberg, I. J.},
     TITLE = {A note on the cyclotomic polynomial},
   JOURNAL = {Mathematika},
    VOLUME = {11},
      YEAR = {1964},
     PAGES = {131--136},
      ISSN = {0025-5793},
        URL = {https://doi.org/10.1112/S0025579300004344},
}

\bib{Ste08}{article}{
    AUTHOR = {Steinberger, John P.},
     TITLE = {Minimal vanishing sums of roots of unity with large
              coefficients},

   JOURNAL = {Proc. Lond. Math. Soc. (3)},
    VOLUME = {97},
      YEAR = {2008},
    NUMBER = {3},
     PAGES = {689--717},
      ISSN = {0024-6115},
       URL = {https://doi.org/10.1112/plms/pdn006},
}

\bib{Wah97}{book}{
    AUTHOR = {Washington, L. C.},
     TITLE = {Introduction to cyclotomic fields},
    SERIES = {Graduate Texts in Mathematics},
    VOLUME = {83},
   EDITION = {Second},
 PUBLISHER = {Springer-Verlag, New York},
      YEAR = {1997},
     PAGES = {xiv+487},
      ISBN = {0-387-94762-0},
}

\end{thebibliography}
\end{document}